\documentclass[12pt]{article} 
\usepackage{latexsym,amssymb,amsmath}
\usepackage{graphicx} 
 \usepackage{wrapfig}
\usepackage{cite}
 \usepackage{url}
\usepackage{fancyhdr}
\usepackage{amsthm}
\usepackage[all]{xy}
\usepackage{fontenc}
\usepackage{listings}
\usepackage[usenames]{color}
\usepackage[top=1in, bottom=1in, left=2.5cm, right=1in]{geometry}

\newtheorem{theorem}{Theorem}[section]
  \newtheorem{proposition}{Proposition}[section]
               \newtheorem{corollary}{Corollary}[section]
               \newtheorem{lemma}{Lemma}[section]
               
               \newtheorem{remark}{Remark}[section]
               
               \def\pf{\par\noindent {\em Proof.}~\par\noindent}
               \def\qed{~\hfill{$\square$}\pagebreak[1]\par\medskip\par}

\newcommand{\R}{{\mathbb R}}

\newcommand{\I}{{\cal I}}

\newcommand{\pD}{{^\psi\!\upa}}
\newcommand{\hD}{{^\varphi\!\upa}}
\newcommand{\wH}{{^{\varphi,\psi}\!\Psi}}
\newcommand{\wPH}{{^{\varphi,\varphi}\!\Psi}}
\newcommand{\upa}{\underline{\partial}}

\newcommand{\cH}{{\cal H}}

\begin{document}
\title{Generalizations of harmonic functions in $\R^m$}
\small{
\author
{Daniel Alfonso Santiesteban; Yudier Pe\~na P\'erez;\\Ricardo Abreu Blaya}
\vskip 1truecm
\date{\small Facultad de Matem\'aticas, Universidad Aut\'onoma de Guerrero, M\'exico.\\ Emails: danielalfonso950105@gmail.com, ypenap88@gmail.com, rabreublaya@yahoo.es}

\maketitle
\begin{abstract}
In recent works, arbitrary structural sets in the non-commutative Clifford analysis context have been used to introduce non-trivial generalizations of harmonic Clifford algebra valued functions in $\R^m$. Being defined as the solutions of elliptic (generally non-strongly elliptic) partial differential equations, $(\varphi,\psi)$-inframonogenic and $(\varphi,\psi)$-harmonic functions do not share the good structure and properties of the harmonic ones. The aim of this paper it to show and clarified the relationship between these classes of functions.     
\end{abstract}

\vspace{0.3cm}

\small{
\noindent
\textbf{Keywords.} Clifford analysis, harmonic functions, structural sets.\\
\noindent
\textbf{Mathematics Subject Classification (2021).} 30G35.}
\section{Introduction}
In recent works \cite{A1,S1} the authors used the so called structural sets to define two subclasses of biharmonic functions in $\R^m$, which may be interpreted as somewhat exotic generalizations of classical harmonic functions.   

A structural set can be though as an orthonormal basis $\psi=\{\psi^1, \psi^2,\dots,\psi^m\}$ of $\R^m$, whose vectors being subjected to the multiplication rules
\[\psi^i\psi^j+\psi^j\psi^i=-2\delta_{ij}\,\,(i,j=1,2,\dots, m)\]
generate the real Clifford algebra $\R_{0,m}$.

The conventional notion of Clifford monogenicity \cite{BDS} refers to the $\R_{0,m}$-valued solutions of the Dirac equation $\upa f=0$, where 
\[
\upa=e_1{\frac{\partial}{\partial x_1}}+e_2{\frac{\partial}{\partial x_2}}+\cdots +e_m{\frac{\partial}{\partial x_m}}
\]
stands for the orthogonal Dirac operator in $\R^m$ constructed with the standard basis $\{e_1,e_2,\dots,e_m\}$. 

More generally, the notion of $\psi$-hyperholomorphic functions \cite{GN1,GN2,SV} ties up the above mentioned monogenicity with arbitrary structural sets. In this way, for a fixed arbitrary structural set $\psi$, an $\R_{0,m}$-valued function $f$ is said to be (left or right respectively) $\psi$-hyperholomorphic if it belongs to $\ker[\pD(\cdot)]$ or $\ker[(\cdot)\pD]$, where
\[
\pD :=\psi^1\frac{\partial}{\partial x_1}+\psi^2\frac{\partial}{\partial x_2}+\cdots +\psi^m\frac{\partial}{\partial x_m}.
\]
As is easily seen, $\pD$ factorizes the Laplace operator in the sense that 
\[\pD^2=\pD\pD=-\Delta.\] 
Having in mind this factorization, the Laplace equation $\Delta f=0$ may be written as $\pD\pD[f]=0$. The consideration of arbitrary and distinct  structural sets $\varphi$, $\psi$, will give us the flexibility to deal with the second order partial differential equation $\hD\pD[f]=0$ whose solutions have been referred to as $(\varphi,\psi)$-harmonic functions (see \cite{S1}). In a more general context we refer to \cite{BDGS}, where higher order equations associated to a finite number of different structural sets has been studied. At the same time the non-commutativity of the product in $\R_{0,m}$ gives rise to the sandwich equation $\hD [f]\pD=0$. The solutions of this last equation represent a natural generalization of the so-called inframonogenic functions studied in \cite{MPS1,MPS2,MAB1,MAB2,MAB3}, which explains the name of $(\varphi,\psi)$-inframonogenic functions adopted for them in \cite{A1}.

In the sequel, when speaking of domains, they will always be assumed to be open and connected sets of $\R^m$. We will use the Greek capital letter $\Omega$ to denote a generic domain. To avoid needless repetition, we specify at once that we use the symbols $\cH(\Omega)$, $\cH_{\varphi,\psi}(\Omega)$ and $\I_{\varphi,\psi}(\Omega)$ to denote the class of harmonic, $(\varphi,\psi)$-harmonic and $(\varphi,\psi)$-inframonogenic functions in $\Omega$, respectively.

Of course, since $\pD^2=\hD^2=-\Delta$, one may conclude that the above three sets are proper subspaces of the space of $\R_{0,m}$-valued biharmonic functions in $\Omega$. The relative positions of these function subspaces are illustrated in  Figure $1$.
\begin{figure}[h]
\centering
\includegraphics[scale=0.06]{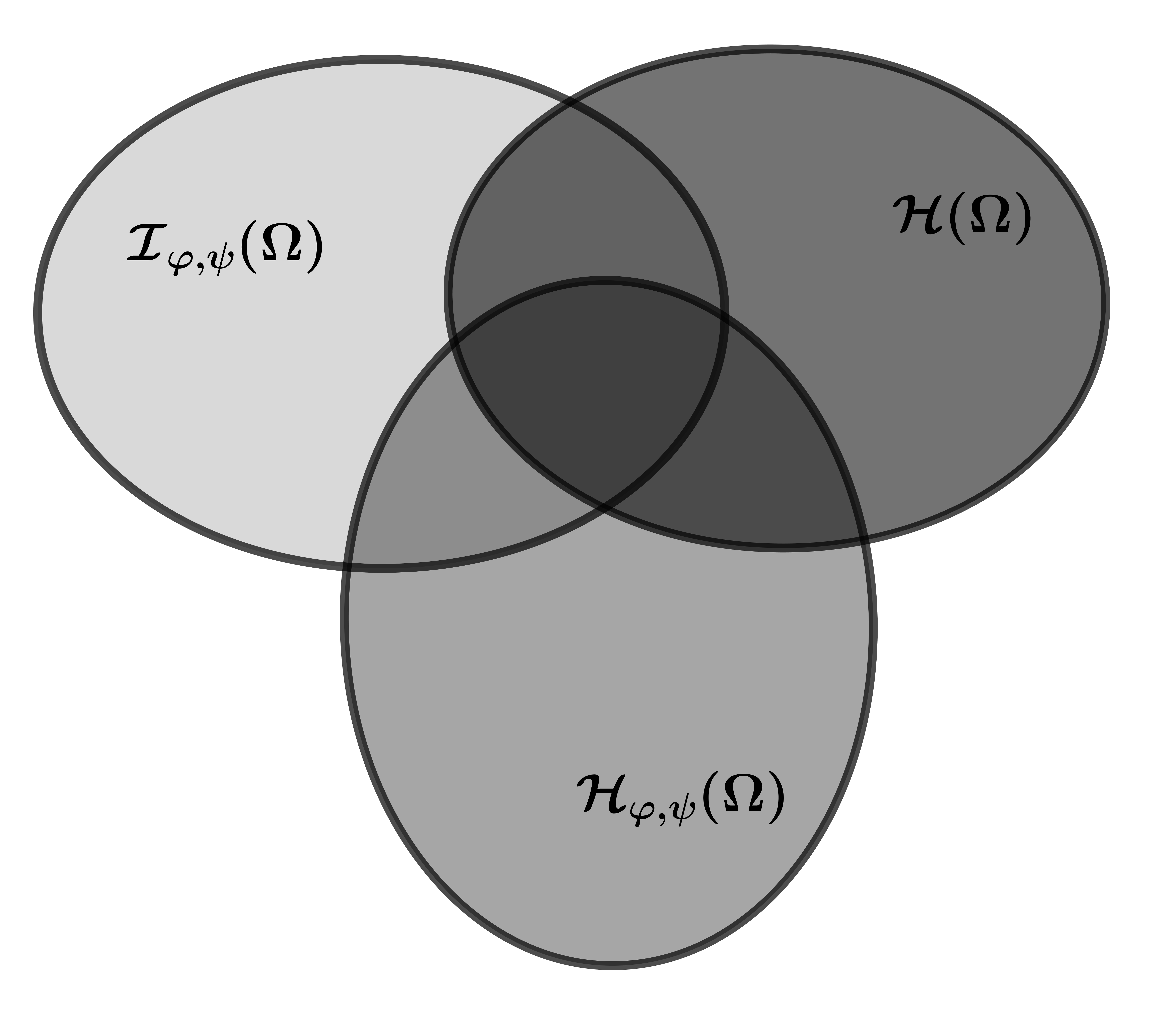}
\caption{The subspaces $\cH(\Omega)$, $\cH_{\varphi,\psi}(\Omega)$ and $\I_{\varphi,\psi}(\Omega)$}
\end{figure}

The following concrete examples help us to give a more complete picture of Figure $1$ when $\Omega=\R^3$ and $\varphi=\{e_1,e_2,e_3\}$, $\psi=\{e_3,e_2,e_1\}$.

The $\R_{0,3}$-valued polynomial $(x_2^2-x_1^2)e_2-2x_1x_2e_3-x_1e_1e_2+x_3e_2e_3$, being (left and right) $\psi$-hyperholomorphic in $\R^3$, obviously belongs to 
\begin{equation}\label{inter}
\cH(\R^3)\cap\cH_{\varphi,\psi}(\R^3)\cap\I_{\varphi,\psi}(\R^3).
\end{equation} 
Direct computations however confirm that the $2$-th degree polynomial $2x_1x_3e_1-x_2e_2-(x_1^2-x_3^2)e_3$ is indeed in the above intersection without being $\psi$-hyperholomorphic in $\R^3$. 

On the other hand, the polynomial $2x_2x_3e_1-(x_1^2+x_2^2)e_2$ belongs to $\cH_{\varphi,\psi}(\R^3)\cap\I_{\varphi,\psi}(\R^3)\setminus\cH(\R^3)$,
meanwhile $x_1x_3e_1+x_2e_2$ and $(x_1x_2+x_2x_3)e_2$ belong to $\cH(\R^3)\cap\I_{\varphi,\psi}(\R^3)\setminus\cH_{\varphi,\psi}(\R^3)$ and $\cH(\R^3)\cap\cH_{\varphi,\psi}(\R^3)\setminus\I_{\varphi,\psi}(\R^3)$, respectively. Similarly, one can easily construct examples of functions which belong to only one of the subspaces appearing in \eqref{inter}.  

In this paper, we introduce a generalization of the operator 
\[\Psi(a)=\sum_{j=1}^m e_j a e_j,\]
which has been extensively used in \cite{ABMM,MMA}. This generalization aids us in finding several relationships between the subspaces $\cH(\Omega)$, $\cH_{\varphi,\psi}(\Omega)$ and $\I_{\varphi,\psi}(\Omega)$. 
\section{Preliminaries and auxiliary results}
We begin by recalling that any Clifford number $a\in\R_{0,m}$ may  be written as $a=\sum_{k=0}^m[a]_k$ where $[a]_k$ is the projection of $a$ onto the subspace of $k$-grade multivectors $\R_{0,m}^{(k)}$ defined by
$$\R_{0,m}^{(k)}=\left\{a\in\R_{0,m}:\;a=\sum_{|A|=k}a_Ae_A,\;\;a_A\in\R\right\}.$$ 
Here $A=\{j_1,\dots,j_k\}\subset\{1,\dots,m\}$ ($j_1<\cdots <j_k$) and $e_A=e_{j_1}\cdots e_{j_k}$.

The real Clifford algebra $\R_{0,m}$ can be split up into two subspaces $\R_{0,m}^+$ and $\R_{0,m}^-$ containing respectively the even and odd multivectors, namely $\R_{0,m}=\R_{0,m}^+\oplus\R_{0,m}^-$. 

In this way, any Clifford number $a$ admits also the unique splitting
\[
a=a_++a_-,\,\,a_{\pm}\in\R_{0,m}^\pm,
\] 
where $a_+$ ($a_-$) are referred to as the even (odd) part of $a$.

We shall use of the following anti-automorphisms in $\R_{0,m}$: the conjugation $a\mapsto\overline{a}$ defined by $\overline{e_i}=-e_i$, and the reversion $a\mapsto\hat{a}$ defined by $\widehat{e_i}=e_i$, $i=1,\dots,m$.

Let $\varphi$ and $\psi$ be two arbitrary structural sets in $\R^m$. For $k=1,\dots,m$, we introduce the operators $\wH_k:\R_{0,m}\to\R_{0,m}$ given by
\[
\wH_k(a):=\sum_{|A|=k}\varphi_A a\widehat{\psi_A},
\]
where $\varphi_A=\varphi^{j_1}\cdots\varphi^{j_k}$, $\psi_A=\psi^{j_1}\cdots\psi^{j_k}$.

For $k=0$ we adopt the convention  $\wH_0:=I$, being the identity operator. 

In particular, if we take $\varphi=\psi=\{e_1,\dots,e_m\}$ and $k=1$ we get
\[
\wH_1(a)=\sum_{j=1}^m e_j a e_j,
\]
which is the above-mentioned operator $\Psi$ introduced in \cite{ABMM}. As proved in \cite{MMA}, when restricting to odd dimension $m$, the operator $\Psi$ becomes a bijection. We will see now that this is a characteristic property of the generalized operators $\wH_1$, a fact that can be regarded a special case of the following general reasoning. 

Let $\{j_1,\dots,j_k\}\subset\{1,\dots,m\}$ and consider the sub-collections $\{\varphi^{j_1},\dots,\varphi^{j_k}\}$ and $\{\psi^{j_1},\dots,\psi^{j_k}\}$ of $\varphi$ and $\psi$, respectively. Now we formulate
\begin{proposition}\label{p1}
If $k$ is odd, then the operator
$$\wH_1^{\{j_1,\dots,j_k\}}:\,a\mapsto\sum_{i=1}^k\varphi^{j_i}a\psi^{j_i},$$
is a real linear bijection on $\R_{0,m}$.
\end{proposition}
\pf We are reduced to prove that 
\[\ker\bigg[\wH_1^{\{j_1,\dots,j_k\}}\bigg]=\{0\}.\]
We proceed by complete induction on $k$. For $k=1$ the conclusion is immediate, since $\varphi^{j_1}$ and $\psi^{j_1}$ are both invertible $\R_{0,m}^{(1)}$-valued constants. 

Let us assume that  
\begin{equation}\label{=0}
\wH_1^{\{j_1,\dots,j_{k+2}\}}(a)=0.
\end{equation}
If we multiply \eqref{=0} on the left by $\varphi^{j_i}$, on the right by $\psi^{j_i}$ and use the identities $\varphi^{j_i}\varphi^{j_i}=-1$, $\psi^{j_i}\psi^{j_i}=-1$, we see that
\[
\left\{\begin{array}{ccc}
a+\varphi^{j_1}\varphi^{j_2}a\psi^{j_2}\psi^{j_1}+\cdots +\varphi^{j_1}\varphi^{j_{k+2}}a\psi^{j_{k+2}}\psi^{j_1}&=&0\\
a+\varphi^{j_2}\varphi^{j_1}a\psi^{j_1}\psi^{j_2}+\cdots +\varphi^{j_2}\varphi^{j_{k+2}}a\psi^{j_{k+2}}\psi^{j_2}&=&0\\
\vdots &\vdots &\vdots\\
a+\varphi^{j_{k+2}}\varphi^{j_1}a\psi^{j_1}\psi^{j_{k+2}}+\cdots +\varphi^{j_{k+2}}\varphi^{j_{k+1}}a\psi^{j_{k+1}}\psi^{j_{k+2}}&=&0.
\end{array}\right.
\]
Subtracting the first equation from the others in the above system, one obtains
\[
\left\{\begin{array}{ccc}
\wH_1^{\{j_3,\dots,j_{k+2}\}}(\varphi^{j_1}a\psi^{j_1})&=&\wH_1^{\{j_3,\dots,j_{k+2}\}}(\varphi^{j_2}a\psi^{j_2})\\
\wH_1^{\{j_2,j_4,\dots,j_{k+2}\}}(\varphi^{j_1}a\psi^{j_1})&=&\wH_1^{\{j_2,j_4,\dots,j_{k+2}\}}(\varphi^{j_3}a\psi^{j_3})\\
\vdots &\vdots &\vdots\\
\wH_1^{\{j_2,\dots,j_{k+1}\}}(\varphi^{j_1}a\psi^{j_1})&=&\wH_1^{\{j_2,\dots,j_{k+1}\}}(\varphi^{j_{k+2}}a\psi^{j_{k+2}}).
\end{array}\right.
\]
Hence, by the induction hypothesis we have
\[\varphi^{j_1}a\psi^{j_1}=\varphi^{j_i}a\psi^{j_i}\;\forall i\in\overline{2,k+2},\]
which together with \eqref{=0} implies that $a=0$.\qed
\begin{corollary}
Let $\varphi$ and $\psi$ be two arbitrary structural sets and suppose that $m$ is odd. Then $\wH_1$ is a bijective $\R$-linear mapping on $\R_{0,m}$. 
\end{corollary}
Some interesting and explicit relations arise when we consider the particular case $\varphi=\psi$. Indeed we have

\begin{proposition}\label{p2}
Let be $a_k\in\R_{0,m}^{(k)}$, then
\begin{equation}\label{wak}
\wPH_j(a_k)=(-1)^{j(k+1)}\sum_{i=\max\{0,j+k-m\}}^{\min\{j,k\}}\binom{m-k}{j-i}\binom{k}{i}(-1)^{i}a_k.
\end{equation}
\end{proposition}
\pf The proof relies on a direct computation and will be omitted.
\begin{remark}\label{R1}
Formula \eqref{wak} generalizes the one derived in \cite{MPS2} (see also \cite[Lema 2.1-(4)]{MAB3}) and it can be rewritten using the well-known Gaussian hypergeometric function ${_2\!F_1}$ as follows
\[
\wPH_j(a_k)=\left\{\begin{array}{cc}
(-1)^{j(k+1)}\binom{m-k}{j}\,{_2\!F_1}(-j,-k,1-j-k+m,-1)a_k,&j+k-m\leq 0,\\
(-1)^{k(j+1)+m}\binom{k}{m-j}\,{_2\!F_1}(j-m,k-m,1+j+k-m,-1)a_k,&j+k-m>0.
\end{array}\right.
\]
\end{remark}
The following relations  are straightforward consequences of the above proposition.
\begin{corollary}\label{c2}
If $m$ is odd, then for a given $a\in\R_{0,m}$ we have
\[
\wPH_j(a)=-\wPH_{m-j}(a),\,\,j=0,\dots,m.
\]
On the contrary if $m$ is even, then we have
\[
\wPH_j([a]_k)=\left\{\begin{array}{ccc}\wPH_{m-j}([a]_k)&\mbox{if $k$ is even},\\
-\wPH_{m-j}([a]_k)&\mbox{if $k$ is odd}.
\end{array}\right.
\]
\end{corollary} 
For further use we introduce the natural notations
\[
\wH_+(a)=\sum_{k-even}\wH_{k}(a),\,\,\wH_-(a)=\sum_{k-odd}\wH_{k}(a).
\]
Returning to Corollary \ref{c2} it is noticed that, when $m$ is odd and $\varphi=\psi$, we have 
\[\wPH_+(a)=-\wPH_-(a).\]
More generally, we have the following
\begin{proposition}\label{p3}
If $m$ is odd, then
\[\wPH_+(a)=-\wPH_-(a)=2^{m-1}([a]_0+[a]_m).\]
On the other hand, if $m$ is even, we have 
\begin{align*}
\wPH_+(a)&=2^{m-1}([a]_0+[a]_m),\\
\wPH_-(a)&=2^{m-1}([a]_m-[a]_0).
\end{align*}
\end{proposition}
\pf The proof consists in the convenient use of the identity 
\[
\varphi^j[\wPH_+(a)]\varphi^j=\wPH_-(a),\,\,j=1,\dots,m.
\]
\qed

\section{Main results}
In this section we state and prove our main results. Generally we will consider functions defined on subsets of $\R^{m}$ and taking values in $\R_{0,m}$. Those functions might be written as $f=\sum_{A} f_A e_A$, where $f_A$ are $\R$-valued functions. The notions of continuity, differentiability and integrability of an $\R_{0,m}$-valued function  have the usual component-wise meaning.

Before going further, one interesting remark is in order. As indicated in \cite[Remark 1]{A1}, the nice property
\[
(f\in\I_{\varphi,\varphi}(\Omega))\Longleftrightarrow ([f]_k\in\I_{\varphi,\varphi}(\Omega),\,k=\overline{0,m})
\]
is no longer true when $\varphi\not=\psi$, and the same disadvantage attaches to $(\varphi,\psi)$-harmonic functions. However, all is not lost, the odd and even parts of $f$ enable us to obtain a weaker version of it.
\begin{proposition}\label{p4}
Let $f=f_++f_-$ be a $C^2$-smooth function in $\Omega$, then  
\[
(f\in\I_{\varphi,\psi}(\Omega))\Longleftrightarrow (f_\pm\in\I_{\varphi,\psi}(\Omega)),
\]
\[
(f\in\cH_{\varphi,\psi}(\Omega))\Longleftrightarrow (f_\pm\in\cH_{\varphi,\psi}(\Omega))
.\]
\end{proposition}
\pf The proof follows after a direct checking of the commutation relations: 
\[
[\hD f\pD]_\pm=\hD[f_\pm]\pD,\,\,[\hD\pD f]_\pm=\hD\pD[f_\pm].
\]\qed 

We continue with a proposition, which involves the Dirac operators $\pD$, $\hD$ and the operator $\wH_1$ introduced in the previous section. We tacitly assume that any required differentiations are well defined.
\begin{proposition}\label{p5}
Let be $f:\R^m\to\R_{0,m}$, then
\begin{itemize}
\item[(i)] $\hD[\wH_1(f)]=-2[f]\pD-\wH_1(\hD[f])$,   $[\wH_1(f)]\pD=-2\hD[f]-\wH_1([f]\pD)$,
\item[(ii)] $\hD[\wH_1(f)]\pD=\wH_1(\hD[f]\pD)$, $\Delta[\wH_1(f)]=\wH_1(\Delta f)$,
\item[(iii)] $\wH_1(\hD\pD[f])=-2\pD[f]\pD-\hD[\wH_1(\pD[f])]$.
\end{itemize}
\end{proposition}
\pf
We shall prove the first statements in $(i)$ and $(ii)$, the second identities may be proved in a quite analogous way. 

\textit{Proof of (i):}
\begin{equation*}
\begin{split}
{\hD}[\wH_1(f)]&=\displaystyle\sum_{\overset{1\leq i,j\leq m}{}}\,\varphi^i\varphi^j(\frac{\partial f}{\partial{x_i}})\psi^j\\
&=\displaystyle\sum_{\overset{i=j}{1\leq i,j\leq m}}\,\varphi^i\varphi^j(\frac{\partial f}{\partial{x_i}})\psi^j+\displaystyle\sum_{\overset{i\not=j}{1\leq i,j\leq m}}\,\varphi^i\varphi^j(\frac{\partial f}{\partial{x_i}})\psi^j\\
&=-\displaystyle\sum_{\overset{i=1}{}}^{m}\,(\frac{\partial f}{\partial{x_i}})\psi^i-\displaystyle\sum_{\overset{i\not=j}{1\leq i,j\leq m}}\,\varphi^j\varphi^i(\frac{\partial f}{\partial{x_i}})\psi^j\\
&=-\displaystyle\sum_{\overset{i=1}{}}^{m}\,(\frac{\partial f}{\partial{x_i}})\psi^i-\left[\displaystyle\sum_{\overset{i,j}{}}\,\varphi^j\varphi^i(\frac{\partial f}{\partial{x_i}})\psi^j+\displaystyle\sum_{\overset{i=1}{}}^{m}\,(\frac{\partial f}{\partial{x_i}})\psi^i\right]\\
&=-2\displaystyle\sum_{\overset{i=1}{}}^{m}\,(\frac{\partial f}{\partial{x_i}})\psi^i-\displaystyle\sum_{\overset{i,j}{}}\,\varphi^j\varphi^i(\frac{\partial f}{\partial{x_i}})\psi^j\\
&=-2[f]{\pD}-\wH_1({\hD}[f]).
\end{split}
\end{equation*}

\textit{Proof of (ii):} By $(i)$ we have
\begin{equation*}
\begin{split}
{\hD}[\wH_1(f)]{\pD}&=-2[f]{\pD}{\pD}-[\wH_1({\hD}[f])]{\pD}\\
&=-2[f]{\pD}{\pD}-\{-2{\hD}{\hD}[f]-\wH_1({\hD}[f]{\pD})\}\\
&=\wH_1({\hD}[f]{\pD}).
\end{split}
\end{equation*}

\textit{Proof of (iii):} The proof is straightforward from $(i)$.
\qed
The following theorem gives a necessary and sufficient condition for a $C^2$-smooth function to be $(\varphi,\psi)$-inframonogenic or $(\varphi,\psi)$-harmonic in a domain of $\R^m$.
\begin{theorem}
Suppose that $m$ is odd and let $f\in C^2(\Omega)$, $\Omega\subset\R^m$. Then
\begin{itemize}
\item[(i)] $f\in\I_{\varphi,\psi}(\Omega)$ if and only if $\wH_1(f)\in\I_{\varphi,\psi}(\Omega)$, 
\item[(ii)] $f\in\cH_{\varphi,\psi}(\Omega)$ if and only if $\pD[f]\pD=\frac{1}{2}\hD[\wH_1(\pD[f])]$ in $\Omega$.
\end{itemize} 
\end{theorem}
\pf The proof follows from the above proposition and the bijectivity of $\wH_1$.\qed

The relationship between $\cH(\Omega)$ and $\I_{\varphi,\psi}(\Omega)$  is clarified by the following
\begin{theorem}\label{main1}
Let $f\in C^2(\Omega)$ be harmonic in $\Omega$. Then $\wH_+(f)$ and $\wH_-(f)$ are both harmonic and $(\varphi,\psi)$-inframonogenic functions in $\Omega$. 
\end{theorem}
This fact can be proved directly from the following remarkable identities.
\begin{lemma}\label{L1}
The following formulas hold
\begin{itemize}
\item[(i)] $\hD\left[\wH_+(f)\right]\pD=\wH_-(\Delta f)$,
\item[(ii)] $\hD\left[\wH_-(f)\right]\pD=\wH_+(\Delta f).$
\end{itemize}
\end{lemma}
\pf By definition it follows that
\begin{align*}
\hD\left[\wH_+(f)\right]\pD=\sum_{k-even}{\hD}[\wH_{k}(f)]{\pD}.
\end{align*}
Moreover, we have
\begin{align}
{\hD}[\wH_k(f)]{\pD}&=\sum_{i=1}^m\varphi^i\dfrac{\partial^2[\wH_k(f)]}{\partial x_i^2}\psi^i+\sum_{\overset{1\leq i,j\leq m}{i\not=j}}\varphi^i\dfrac{\partial^2 [\wH_{k}(f)]}{\partial x_i\partial x_j}\psi^j\label{a+b}.
\end{align}
If $\wH_{k}(f)$ is expanded in the second summand of \eqref{a+b}, then  we obtain
\begin{align*}\label{reltrian}
\sum_{\overset{1\leq i,j\leq m}{i\not=j}}\varphi^i\dfrac{\partial^2[\wH_k(f)]}{\partial x_i\partial x_j}\psi^j&=\sum_{\overset{1\leq i,j\leq m}{i\not=j}}\sum_{|A|=k}\varphi^i\varphi_A\dfrac{\partial^2f}{\partial x_i\partial x_j}\widehat{\psi}_A\psi^j.
\end{align*}
Hence
\[
\sum_{k-even}\sum_{\overset{1\leq i,j\leq m}{i\not=j}}\varphi^i\dfrac{\partial^2[\wH_{k}(f)]}{\partial x_i\partial x_j}\psi^j=\sum_{k-even}\sum_{\overset{1\leq i,j\leq m}{i\not=j}}\sum_{|A|=k}\varphi^i\varphi_A\dfrac{\partial^2f}{\partial x_i\partial x_j}\widehat{\psi}_A\psi^j=0,
\]
where the last cancellation may be proved directly by laborious computation.

Now we turn to the first summand in \eqref{a+b}. We have
\[
\sum_{i=1}^m\varphi^i\dfrac{\partial^2[\wH_k(f)]}{\partial x_i^2}\psi^i=\sum_{i=1}^m\sum_{|A|=k}\varphi^i\varphi_A\dfrac{\partial^2 f}{\partial x_i^2}\widehat{\psi}_A\psi^i.
\]
Observe that for all $i=1,\dots,m$:
\begin{align*}
\sum_{|A|=k}\varphi^i\varphi_A\dfrac{\partial^2 f}{\partial x_i^2}\widehat{\psi}_A\psi^i&=\sum_{\overset{|A|=k}{A\ni i}}\varphi^i\varphi_A\dfrac{\partial^2 f}{\partial x_i^2}\widehat{\psi}_A\psi^i+
\sum_{\overset{|A|=k}{A\not\ni i }}\varphi^i\varphi_A\dfrac{\partial^2 f}{\partial x_i^2}\widehat{\psi}_A\psi^i\\
&=\sum_{\overset{|A|=k}{A\ni i}}\varphi_{A\setminus\{i\}}\dfrac{\partial^2 f}{\partial x_i^2}\widehat{\psi}_{A\setminus\{i\}}+
\sum_{\overset{|A|=k}{A\not\ni i }}\varphi_{A\cup\{i\}}\dfrac{\partial^2 f}{\partial x_i^2}\widehat{\psi}_{A\cup\{i\}}.
\end{align*}
On the other hand, for $k\ge 2$ 
\begin{align*}
\sum_{|A|=k}&\varphi^i\varphi_A\dfrac{\partial^2 f}{\partial x_i^2}\widehat{\psi}_A\psi^i+\sum_{|A|=k-2}\varphi^i\varphi_A\dfrac{\partial^2 f}{\partial x_i^2}\widehat{\psi}_A\psi^i\\&=\sum_{\overset{|A|=k}{A\ni i}}\varphi_{A\setminus\{i\}}\dfrac{\partial^2 f}{\partial x_i^2}\widehat{\psi}_{A\setminus\{i\}}+\sum_{\overset{|A|=k}{A\not\ni i }}\varphi_{A\cup\{i\}}\dfrac{\partial^2 f}{\partial x_i^2}\widehat{\psi}_{A\cup\{i\}}\\&+\sum_{\overset{|A|=k-2}{A\ni i}}\varphi_{A\setminus\{i\}}\dfrac{\partial^2 f}{\partial x_i^2}\widehat{\psi}_{A\setminus\{i\}}+
\sum_{\overset{|A|=k-2}{A\not\ni i }}\varphi_{A\cup\{i\}}\dfrac{\partial^2 f}{\partial x_i^2}\widehat{\psi}_{A\cup\{i\}}\\&=
\sum_{\overset{|A|=k-2}{A\ni i}}\varphi_{A\setminus\{i\}}\dfrac{\partial^2 f}{\partial x_i^2}\widehat{\psi}_{A\setminus\{i\}}+\sum_{|A|=k-1}\varphi_A\dfrac{\partial^2 f}{\partial x_i^2}\widehat{\psi}_A+\sum_{\overset{|A|=k}{A\not\ni i }}\varphi_{A\cup\{i\}}\dfrac{\partial^2 f}{\partial x_i^2}\widehat{\psi}_{A\cup\{i\}},
\end{align*}
where use has been made of the obvious equality
\[\sum_{\overset{|A|=k}{A\ni i}}\varphi_{A\setminus\{i\}}\dfrac{\partial^2 f}{\partial x_i^2}\widehat{\psi}_{A\setminus\{i\}}+\sum_{\overset{|A|=k-2}{A\not\ni i }}\varphi_{A\cup\{i\}}\dfrac{\partial^2 f}{\partial x_i^2}\widehat{\psi}_{A\cup\{i\}}=\sum_{|A|=k-1}\varphi_A\dfrac{\partial^2 f}{\partial x_i^2}\widehat{\psi}_A.
\]
A repeated use of the above argument leads to
\[
\sum_{k-even}\sum_{|A|=k}\varphi^i\varphi_A\dfrac{\partial^2 f}{\partial x_i^2}\widehat{\psi}_A\psi^i=\sum_{\kappa-odd}\sum_{|A|=\kappa}\varphi_A\dfrac{\partial^2 f}{\partial x_i^2}\widehat{\psi}_A.
\]
Therefore
\begin{align*}
\hD\left[\wH_+(f)\right]\pD &=\sum_{i=1}^m\sum_{k-even}\sum_{|A|=k}\varphi^i\varphi_A\dfrac{\partial^2 f}{\partial x_i^2}\widehat{\psi}_A\psi^i=\sum_{\kappa-odd}\sum_{|A|=\kappa}\sum_{i=1}^m\varphi_A\dfrac{\partial^2 f}{\partial x_i^2}\widehat{\psi}_A\\&=\sum_{\kappa-odd}\sum_{|A|=\kappa}\varphi_A[\Delta f]\widehat{\psi}_A=\wH_-(\Delta f).
\end{align*}
The second identity may be proved in a quite analogous way.\qed 
{\emph{Proof of Theorem \ref{main1}}.} If $f$ is harmonic, then the harmonicity of $\wH_+(f)$ follows directly from the obvious identity $\Delta[\wH_+(f)]=\wH_+(\Delta f)$. That the function $\wH_+(f)$ is also $(\varphi,\psi)$-inframonogenic is clear from Lemma \ref{L1}-$(i)$. In a similar fashion, the function $\wH_-(f)$ is dealt with.\qed 

As we shall now see, the conclusion of Theorem \ref{main1} remains true if $\lq\lq$harmonic" is replaced by $\lq\lq(\varphi,\psi)$-inframonogenic". More precisely:
\begin{theorem}\label{main2}
Let $f$ be $(\varphi,\psi)$-inframonogenic in $\Omega$. Then $\wH_+(f)$ and $\wH_-(f)$ are both harmonic and $(\varphi,\psi)$-inframonogenic functions in $\Omega$. 
\end{theorem}  
\pf The proof is obtained as a combination of Proposition \ref{p5}-$(ii)$, Lemma \ref{L1} and the following recursive rule 
\begin{equation}\label{recomega}
(m-k+1)\wH_{k-1}(f)+(k+1)\wH_{k+1}(f)=\wH_1(\wH_k(f)),\;\;k=1,\dots,m-1,
\end{equation}
which may be verified by direct checking.

Indeed, if $f\in\I_{\varphi,\psi}(\Omega)$ then by Proposition \ref{p5}-$(ii)$ one has $\wH_1(f)\in\I_{\varphi,\psi}(\Omega)$. This last fact and \eqref{recomega} together say that $\wH_2(f)\in\I_{\varphi,\psi}(\Omega)$. Thus, this argument may be applied repeatedly to obtain that 
$\wH_k(f)$ ($k=1,\dots,m$) and so $\wH_\pm(f)$ are $(\phi,\psi)$-inframonogenic functions in $\Omega$ as well. With this in hand, the harmonicity of $\wH_\pm(f)$ is a straightforward consequence of Lemma \ref{L1}.\qed

Summarizing, we have proved that the operators $\wH_\pm$ define a pair of linear mappings from $\cH(\Omega)\cup\I_{\varphi,\psi}(\Omega)$ to $\cH(\Omega)\cap\I_{\varphi,\psi}(\Omega)$ (See Figure $2$). 
\begin{figure}[h]
\centering
\includegraphics[scale=0.1]{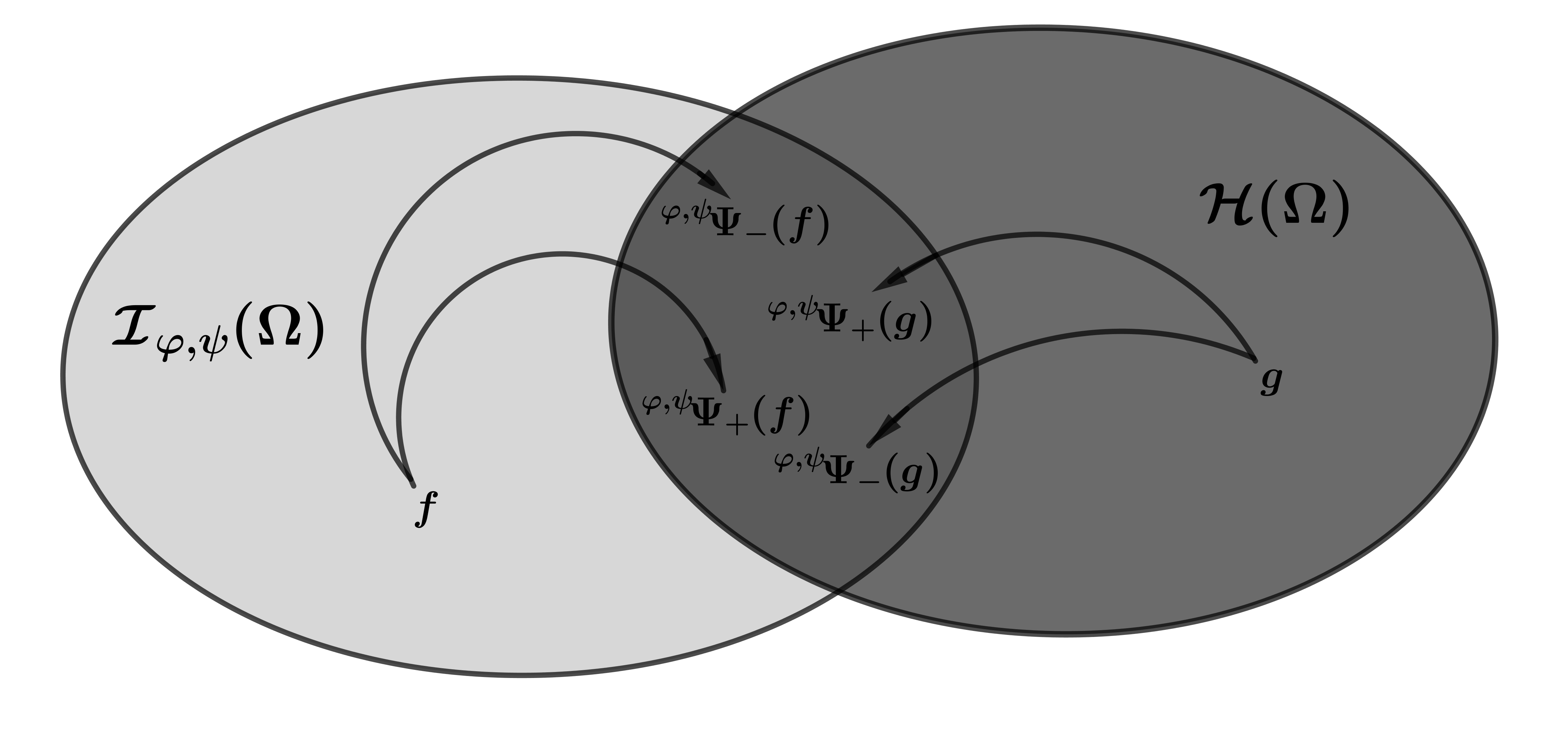}
\caption{The operators $\wH_\pm$ as mappings from $\cH(\Omega)\cup\I_{\varphi,\psi}(\Omega)$ to $\cH(\Omega)\cap\I_{\varphi,\psi}(\Omega)$}
\end{figure}

It is worthwhile to see how our main theorems work in the $2$-dimensional case.  Let $\varphi$ and $\psi$ be two structural sets in $\R_{0,2}$ and consider a given function $f=f_0+f_1\psi_1+f_2\psi_2+f_{12}\psi_1\psi_2$.  The transition matrices from the basis $\psi$ to the basis $\varphi$ are of the form 
\begin{equation}\label{t1}
\left(\begin{array}{cccc}
c_{1}&-c_{2}\\c_{2}&c_{1}
\end{array}\right)
\end{equation} 
or
\begin{equation}\label{t2}
\left(\begin{array}{cccc}
c_{1}&c_{2}\\c_{2}&-c_{1}
\end{array}\right),
\end{equation} 
where $c_1^2+c_2^2=1$. 

In the first case \eqref{t1}, we have
\[
\wH_+(f)=2f_0+2f_{12}\psi^1\psi^2,\,\,\wH_-(f)=-2(c_1f_0+c_2f_{12})+2(c_1f_{12}-c_2f_0)\psi^1\psi^2.
\]
So in this case Theorems \ref{main1} and \ref{main2} say that if $f=f_++f_-$ is harmonic or $(\varphi,\psi)$-inframonogenic, then $f_+$ is simultaneously harmonic and $(\varphi,\psi)$-inframonogenic.

On the other hand, in case \eqref{t2} one has 
\[
\wH_+(f)=2f_1\psi^1+2f_2\psi^2,\,\,\wH_-(f)=-2(c_1f_1+c_2f_2)\psi^1+2(c_1f_2-c_2f_1)\psi^2.
\]
Therefore, in this situation, our theorems ensure that if $f=f_++f_-$ is harmonic or $(\varphi,\psi)$-inframonogenic, then $f_-$ is simultaneously harmonic and $(\varphi,\psi)$-inframonogenic.

Of course, the two-dimensional case is very special. In general such a nice interpretation of our theorems is no longer available.

Finally we mention that the converse of the above theorems is not true. As is easily seen from Proposition \ref{p3}, a general $C^2$-smooth function $f$ (neither harmonic nor $(\varphi,\varphi)$-inframonogenic) can be constructed in such a way that $[f]_0$ and $[f]_m$, hence $\wPH_+(f)$, be harmonic and $(\varphi,\varphi)$-inframonogenic, simultaneously. 

\end{document}